\documentclass{article}
\usepackage{biblatex}
\usepackage{mathrsfs}
\usepackage{graphicx}

\title{A New Proof for The Transformation Laws Of Jacobi Theta Functions}
\author{Maher Me'meh $\&$ Ali Saraeb }
\date{}

\begin{document}

\maketitle
\begin{flushleft}
\textbf{Abstract.} We prove the transformation laws of the four Jacobi theta functions using Gordon's proof for the transformation law of the Dedekind eta function.\end{flushleft}

\begin{center}
    1.~ INTRODUCTION
\end{center}
The Jacobi theta functions form an important class of functions in elliptic function theory. These are quasi-periodic entire functions that have been originally formulated by Jacobi and extensively studied by many mathematicians in different disciplines . These functions satisfy numerous identities showing up in the fields of ODEs, abelian varieties, moduli spaces, quadratic forms and quantum field theory.
\\
We define the Jacobi theta functions adopting the notation as used in [2], $\tau \in ~\mathbf{H}$ and $z \in \mathbf{C}$ where $q=e^{\pi i \tau}$.
\[\theta_1(z,\tau)=-i\sum_{n=-\infty}^{\infty} (-1)^n q^{(n+1/2)^2}e^{(2n+1)i\pi z}\]
\[\theta_2(z,\tau)=\sum_{n=-\infty}^{\infty} q^{(n+1/2)^2}e^{(2n+1)\pi i z}\]
\[\theta_3(z,\tau)=\sum_{n=-\infty}^{\infty} q^{n^2}e^{2n\pi i z}\]
\[\theta_4(z,\tau)=\sum_{n=-\infty}^{\infty} (-1)^n q^{n^2}e^{2pi i n z}\]
\\
Many proofs have been provided for the transformation laws of these functions, for instance  using Seigel's method for the Dedekind eta function $\eta(\tau)$ in [4] to prove the inversion formula of $\theta_3$. Also some proofs have been provided for theta functions of higher degrees using the theory of quadratic forms, see [5].
\vspace{0.2cm}
\\
In this paper, we present a new proof of the transformation law of $\theta_1$ under $\Gamma$, and $\theta_2,~\theta_3,~\theta_4$ under $\Gamma(2)$. Our proof is inspired by Basil Gordon's proof of the Dedekind eta function $\eta(\tau)$ in [1]. We intend to establish newer and lighter proofs to already existing standard results of these $\theta$-functions.\\
The idea of the proof is the following: If the transformation law is satisfied for at least one matrix A in $\Gamma$  and works for $AS$ and $AT^m$ where $S$ and $T^m$ are the generators of $\Gamma$  then it works for all matrices in the modular group $\Gamma$. The same proof would follow for the cases of $\theta_2,~\theta_3,~\theta_4$ under the generators of $\Gamma(2).$\\
In this paper we give the full proof for $\theta_1$ and $\theta_3$ only. Since $\theta_1,\theta_2,\theta_3$ and $\theta_4$ are related to each other through relations described in [3-X], the same proof follows for $\theta_2$ and $\theta_4$. 
\vspace{0.3cm}
\\
\begin{center}
    2. THE TRANSFORMATION LAW FOR $\theta_1$
\end{center}
We now determine the transformation law for $\theta_1$ on the full modular group $\Gamma$. The transformation law of $\theta_1$ for the matrices $T^m=\left(\begin{array}{cc}
   1  &  m\\
   0  &1 \\
\end{array}\right)$ $\&$ $S=\left(\begin{array}{cc}
   0  &  -1\\
   1  &0 \\
\end{array}\right)$, see [3-X] are given by:
\[~~~~~~~~\theta_1(z,\tau +m)=e^{\pi i m/4}\theta_1(z,\tau)~~~~~~~~~(1)\]
\[~~~~~~\theta_1(\frac{z}{\tau},\frac{-1}{\tau})=-i\sqrt{-i\tau}e^{\frac{\pi i z^2}{\tau}}\theta_1(z,\tau).~~~~~~(2)\]
The proof will follow as such: If the transformation law were to hold for one matrix $A\in \Gamma$ and it holds for $AT^m$ and $AS$, where T and S are described above, then it holds for all matrices in $\Gamma$. The key observation is that it holds for at least one $A$, namely $S=\left(\begin{array}{cc}
   0  &  -1\\
   1  &0 \\
\end{array}\right)$, where $c>0$.
\vspace{0.35cm}
\\
\textbf{Theorem 1.} If $A=\left(\begin{array}{cc}
   a  &  b\\
   c  &d \\
\end{array}\right)~\in~\Gamma$ such that $c>0$ then the transformation of $\theta_1$ is given by
\[~~~~~~~~~~~\theta_1\left(\frac{z}{c\tau+d},\frac{a\tau+b}{c\tau+d}\right)=\epsilon_1(A)\left(-i(c\tau+d)\right)^{1/2}e^{\frac{\pi i c z^2}{c\tau+d}}\theta_1(z,\tau)~~~~~~~~~~~(3)\]
where $\displaystyle \epsilon_1(A)=-i\epsilon^3=\left\lbrace \begin{array}{ll}
\left(\frac{d}{c}\right)i^{(c-3)/2}e^{(\pi i/4)c(a+d)}~~~~c~odd \\
\left(\frac{c}{d}\right)e^{\pi i/4}i^{(1-d)/2} e^{(\pi i /4)d(b-c)}~~~~d~odd
\end{array}\right.$\\
(check [1-X]).
\vspace{0.3cm}
\\
Here $\epsilon$ appears in the transformation law of the Dedekind eta function, see [3-X].
\[\eta(A\tau)=\epsilon(A) (-i(c\tau+d))^{1/2}\eta(\tau)\]
where
\[\epsilon(A)=exp(\pi i (\frac{a+d}{12c}-s(d,c))),\]
\[\displaystyle s(h,k)=\sum_{r=1}^{k-1}\frac{r}{k}\left(\frac{hr}{k}-\left[\frac{hr}{k}\right]-\frac{1}{2}\right)\] is the Dedekind sum for $k>0$ and $(k,h)=1$. \\
Hence \[\epsilon_1(A)=-i\epsilon^3=-i.exp\left(3\pi i\left(\frac{a+d}{12c}-s(d,c)\right)\right).\]
Now in order to prove Theorem 1, we need the following two Lemmas.

\vspace{0.3cm}
\textbf{Lemma 1.} $\epsilon_1(AT^m)=\epsilon_1(A).e^{\frac{\pi i m}{4}}$
\\
\textit{Proof.} We have \[AT^m=\left(\begin{array}{cc}
   a  &  b\\
   c  &d \\
\end{array}\right).\left(\begin{array}{cc}
   1  &  m\\
   0  &1 \\
\end{array}\right)=\left(\begin{array}{cc}
   a  & am+b \\
   c  &cm+d \\
\end{array}\right).\]
As a result, \[\epsilon(AT^m)=-i exp(3\pi i (\frac{a+cm+d}{12c}-s(cm+d,c))).\] Using the well-known property of the Dedekind sum, $s(cm+d,c)=s(d,c)$, we get
\[\epsilon_1(AT^m)=-i exp(3\pi i (\frac{a+cm+d}{12c}-s(cm+d,c)))\]
\[~~~~=-i exp(3\pi i (\frac{a+d}{12c}-s(d,c)+\frac{m}{12}))\]
\[=-i\epsilon(A).e^{\pi i m/4}=\epsilon_1(A).e^{\pi i m/4}.\]
\vspace{0.3cm}
\\
\textbf{Lemma 2.} $\epsilon_1(AS)=\left\lbrace\begin{array}{ll}
e^{-3\pi i/4}\epsilon_1(A)~~~~~~~if~d>0\\
e^{3\pi i/4}\epsilon_1(A)~~~~~~~if~ d<0\\
\end{array}\right.$\\

\textit{Proof.} First we treat the case when $d>0$, we use $S=\left(\begin{array}{cc}
   0  &  -1\\
   1  &0 \\
\end{array}\right)$, so we have $AS=\left(\begin{array}{cc}
   b  &  -a\\
   d  &-c \\
\end{array}\right)$ and \[~~~~~~~~~~~~\epsilon_1(AS)=i exp\left(3\pi i \left(\frac{b-c}{12d}-s(-c,d)\right)\right)=i exp\left(3\pi i \left(\frac{b-c}{12d}+s(c,d)\right)\right)~~~~~~~~~~~(*)\] using the property $s(-h,k)=-s(h,k)$.\\
We now use the reciprocity law of the Dedekind sum
\[s(d,c)+s(c,d)=\frac{c}{12d}+\frac{d}{12c}-\frac{1}{4}+\frac{1}{12cd}.\]
Replacing $1=ad-bc$, we have
\[s(d,c)+s(c,d)=\frac{c}{12d}+\frac{d}{12c}-\frac{1}{4}+\frac{ad-bc}{12cd}.\]
we obtain
\[s(c,d)=-s(d,c)+\frac{c}{12d}+\frac{d}{12c}-\frac{1}{4}+\frac{a}{12c}-\frac{b}{12d}\]
\[~~~~~~=-s(d,c)+\frac{c-b}{12d}+\frac{a+d}{12c}-\frac{1}{4}.~~~~~~~~~~~\]
Replacing in $\epsilon_1(AS)$ in (*), we have:
\[\epsilon_1(AS)=i exp(3\pi i (\frac{b-c}{12d}-s(d,c)+\frac{c-b}{12d}+\frac{a+d}{12c}-\frac{1}{4}))~\]
\[~~~~~~~~~=i exp(3\pi i (\frac{a+d}{12c}-s(d,c))).e^{-3\pi i/4}~~~~~~~~~~~~~~~~~~\]\[~~~~~~~~~=\epsilon_1(A).e^{-3\pi i/4}.~~~~~~~~~~~~~~~~~~~~~~~~~~~~~~~~~~~~~~~~~~~~\]
Now for the case when $d<0$ we use $S=\left(\begin{array}{cc}
  0   &1  \\
-1     & 0\\
\end{array}\right)$ and so $AS=\left(\begin{array}{cc}
  -b   & a  \\
-d     & c\\
\end{array}\right)$ so that $-d>0$. Thus we have
\[\epsilon_1(AS)=i. exp(3\pi i (\frac{b-c}{12d}-s(c,-d))).\]
Using again the reciprocity law, we get
\[s(c,-d)+s(-d,c)=-\frac{c}{12d}-\frac{d}{12c}-\frac{1}{4}-\frac{ad-bc}{12cd}.\]
Hence
\[s(c,-d)=-s(-d,c)-\frac{c}{12d}-\frac{d}{12c}-\frac{1}{4}-\frac{a}{12c}+\frac{b}{12d}.\]
Substituting again, we obtain
\[\epsilon_1(AS)=i. exp(3\pi i(\frac{b-c}{12d}+s(-d,c)+\frac{c-b}{12d}+\frac{a+d}{12c}+\frac{1}{4}))),\]
\[~~~~~~~~=i. exp(3\pi i(\frac{a+d}{12c}-s(d,c)+\frac{1}{4})))=\epsilon_1(A)e^{3\pi i/4}~~~~~\]
as desired.
\vspace{0.25cm}
\\
We now present the proof of our first thoerem.\\

\textit{Proof of Theorem 1.} 
We follow Gordon's proof as in [1]:\\
Substitute $\tau\to T^m\tau=\tau+m$ in (3) and we have
\[\theta_1\left(\frac{z}{cT^m\tau+d},\frac{aT^m\tau+b}{cT^m\tau+d}\right)=\epsilon_1(A)(-i(cT^m\tau+d))^{1/2}e^{\frac{\pi i c z^2}{cT^m\tau+d}}\theta_1(z,T^m\tau).\]
We get
\[\theta_1\left(\frac{z}{c\tau+cm+d},\frac{a\tau+am+b}{c\tau+cm+d}\right)=\epsilon_1(A)(-i(c\tau+cm+d))^{1/2}e^{\frac{\pi i c z^2}{c\tau+cm+d}}\theta_1(z,\tau+m).~~~~(4)\]
We want to show that (4) is equivalent to (3) when $AT^m=\left(\begin{array}{cc}
   a  & am+b \\
   c  &cm+d \\
\end{array}\right)$ is applied, i.e
\[\theta_1\left(\frac{z}{c\tau+cm+d},\frac{a\tau+am+b}{c\tau+cm+d}\right)=\epsilon_1(AT^m)(-i(c\tau+cm+d))^{1/2}e^{\frac{\pi i c z^2}{c\tau+cm+d}}\theta_1(z,\tau)\]
which is what we have if we use (1) when replacing $\theta_1(z,\tau+m)=e^{\pi i m/4}\theta_1(z,\tau)$ in (4).\\
As a result,
\[\theta_1\left(\frac{z}{c\tau+cm+d},\frac{a\tau+am+b}{c\tau+cm+d}\right)=\epsilon_1(A)(-i(c\tau+cm+d))^{1/2}e^{\frac{\pi i c z^2}{c\tau+cm+d}}.e^{\pi i m/4}\theta_1(z,\tau)\]
\[~~~~~~~~~~~~~~~~~~~~~~~~~~~~~~~~~~~~~~~~~~~~~~~~~~~~~~~~~~~~~~~=\epsilon_1(AT^m)(-i(c\tau+cm+d))^{1/2}e^{\frac{\pi i c z^2}{c\tau+cm+d}}\theta_1(z,\tau)\] using Lemma 1.
\vspace{0.2cm}\\
For the case of S, we treat first the case when $d>0$ and we substitute $\tau \to S\tau=-\frac{1}{\tau}$ in (3) to prove that it is equivalent to applying the matrix $AS=\left(\begin{array}{cc}
   b  &  -a\\
   d  &-c \\
\end{array}\right)$ which will give us
\[~~~~~~~\theta_1\left(\frac{z}{d\tau-c},\frac{b\tau-a}{d\tau-c}\right)=\epsilon_1(AS)(-i(d\tau-c))^{1/2}e^{\frac{\pi i d z^2}{d\tau-c}}\theta_1(z,\tau).~~~~~~~~~~~~(5)\]
Thus,
\[\theta_1\left(\frac{z}{c(-\frac{1}{\tau})+d},\frac{b\tau-a}{d\tau-c}\right)=\epsilon_1(A)\left(-i(c\left(-\frac{1}{\tau}\right)+d)\right)^{1/2}e^{\frac{\pi i c z^2}{c(-\frac{1}{\tau})+d}}\theta_1(z,-\frac{1}{\tau}).\]
As a result,
\[\theta_1\left(\frac{z\tau}{d\tau-c},\frac{b\tau-a}{d\tau-c}\right)=\epsilon_1(A)\left(\frac{-i}{\tau}(d\tau-c)\right)^{1/2}e^{\frac{\pi i c z^2\tau}{d\tau-c}}\theta_1(z,\frac{-1}{\tau})\]
which is not exactly (5).To restore back the same lattices, we do the change of variable $z\to \frac{z}{\tau}$ to get
\[\theta_1\left(\frac{z}{d\tau-c},\frac{b\tau-a}{d\tau-c}\right)=\epsilon_1(A)\left(\frac{-i}{\tau}(d\tau-c)\right)^{1/2}e^{\frac{\pi i c z^2}{\tau(d\tau-c)}}\theta_1(\frac{z}{\tau},\frac{-1}{\tau}).\]
Now we use (2) where \[\theta_1(\frac{z}{\tau},\frac{-1}{\tau})=-i\sqrt{i\tau}e^{\frac{\pi i z^2}{\tau}}\theta_1(z,\tau)\] to get
\[\theta_1\left(\frac{z}{d\tau-c},\frac{b\tau-a}{d\tau-c}\right)=\epsilon_1(A)\left(\frac{-i}{\tau}(d\tau-c)\right)^{1/2}e^{\frac{\pi i c z^2}{\tau(d\tau-c)}}. -i(-i\tau)^{1/2}e^{\frac{\pi i z^2}{\tau}}\theta_1(z,\tau)\]
\[~~~~~~~~~~~~~~=\epsilon_1(A)e^{-3\pi i/4}\left(-i(d\tau-c)\right)e^{\frac{\pi i z^2}{\tau}(\frac{c+d\tau-c}{d\tau-c})}\theta_1(z,\tau)\]
\[=\epsilon(AS)\left(-i(d\tau-c)\right)e^{\frac{\pi iz^2d}{d\tau-c}}\theta_1(z,\tau)~~~\]
using Lemma 2, leading us to (5).\\
For the case when $d<0$, we do a slight change taking again $S=\left(\begin{array}{cc}
  0   &1  \\
-1     & 0\\
\end{array}\right)$. We have $AS=\left(\begin{array}{cc}
  -b   & a  \\
-d     & c\\
\end{array}\right)$ to ensure $-d>0$.
Again imitating the same proof as above we want to obtain
\[ \theta_1\left(\frac{z}{-d\tau+c},\frac{-b\tau+a}{-d\tau+c}\right)=\epsilon_1(AS)\left(-i(-d\tau+c)\right)^{1/2}e^{\frac{-\pi i d  z^2}{-d\tau+c}}\theta_1(z,\tau).~~~~~~~~~~(6)\]
Substituting again $\tau\to -1/\tau$ in (3), we get
\[\theta_1\left(\frac{-z\tau}{-d\tau+c},\frac{-b\tau+a}{-d\tau+c}\right)=\epsilon_1(A)\left(\frac{-i}{-\tau}(-d\tau+c)\right)^{1/2}e^{\frac{-\pi i c z^2\tau}{-d\tau+c}}\theta_1(z,\frac{-1}{\tau}).\]
Doing again the change of variable $z\to -z/\tau$ we get
\[\theta_1\left(\frac{z}{-d\tau+c},\frac{-b\tau+a}{-d\tau+c}\right)=\epsilon_1(A)\left(\frac{-i}{-\tau}(-d\tau+c)\right)^{1/2}e^{\frac{-\pi i c z^2}{\tau(-d\tau+c)}}\theta_1\left(\frac{-z}{\tau},\frac{-1}{\tau}\right).~~~~~~~~~~~(7)\]
Now note that $\theta_1$ is an odd function in terms of z since we can write $\theta_1$ as (see [3-X]):
\[\theta_1(z,\tau)=2\sum_{m=0}^{\infty}(-1)^m q^{(m+1/2)^2}sin\left((2m+1)\pi z\right)\]
so $\theta_1(-z,\tau)=-\theta_1(z,\tau)$. Hence \[~\theta_1(\frac{-z}{\tau},\frac{-1}{\tau})=i\sqrt{-i\tau}e^{\frac{\pi i z^2}{\tau}}\theta_1(z,\tau).\] Substituting in (7), we obtain
\[\theta_1\left(\frac{z}{-d\tau+c},\frac{-b\tau+a}{-d\tau+c}\right)=\epsilon_1(A)\left(\frac{-i}{-\tau}(-d\tau+c)\right)^{1/2}e^{\frac{-\pi i c z^2}{\tau(-d\tau+c)}}.i(-i\tau)^{1/2}e^{\frac{\pi iz^2}{\tau}}\theta_1(z,\tau)\]
\[~~~~=\epsilon_1(A)e^{3\pi i/4}\left(-i(-d\tau+c)\right)^{1/2}e^{\frac{-\pi i z^2}{-d\tau+c}}\theta_1(z,\tau)\]
leading to (6).\\
Since every matrix in $\Gamma$ can be expressed as $A=T^{n_1}ST^{n_2}S...ST^{n_k}$, but also $T=ST^{-1}ST^{-1}S$, then every $A$ can be expressed as $ST^{m_1}ST^{m_2}...ST^{m_r}$ (see [1-III]).
And since it has been proven for $S=\left(\begin{array}{cc}
   0  &  -1\\
   1  &0 \\
\end{array}\right)$, it follows, from above, that the functional equation (3) holds for every $A\in \Gamma$ with $c>0$, as desired.
\vspace{1.5cm}
\\
\begin{center}
    3.~ FOR THE REMAINING THETA FUNCTIONS
\end{center}

As it has been shown $\theta_1$ transforms into itself under elements in $\Gamma$, however this is not the case for $\theta_{2}$, $\theta_3$, $\theta_4$ and that's why we look into the transformation under elements of $\Gamma(2)$\\
\textbf{Theorem 2.} Let $A=\left(\begin{array}{cc}
   a  &  b\\
   c  &d \\
\end{array}\right)\equiv \left(\begin{array}{cc}
   1  & 0 \\
   0  & 1 \\
\end{array}\right) (mod~2) \in \Gamma(2)$ with $c>0$, then, (see [3-6.6]):
\[\theta_2\left(\frac{z}{c\tau+d},\frac{a\tau+b}{c\tau+d}\right)=i^{(d-1)(c/2-1)+c/2}\epsilon_1(A)(c\tau+d)^{1/2}e^{\frac{\pi i z^2 c}{c\tau+d}}\theta_2(z,\tau)~~~~~~~~~~~~\]
\[\theta_3\left(\frac{z}{c\tau+d},\frac{a\tau+b}{c\tau+d}\right)=i^{(d-1)(c/2+1)-\frac{b}{2}(a)+c/2}\epsilon_1(A)(c\tau+d)^{1/2}e^{\frac{\pi i z^2 c}{c\tau+d}}\theta_3(z,\tau)~~~~~(8)\]
\[\theta_4\left(\frac{z}{c\tau+d},\frac{a\tau+b}{c\tau+d}\right)=i^{(a-1)(b/2-1)-b/2}\epsilon_1(A)(c\tau+d)^{1/2}e^{\frac{\pi i z^2 c}{c\tau+d}}\theta_4(z,\tau)~~~~~~~~~~~~~\]
\vspace{0.3cm}
\\
where $\epsilon_1(A)=e^{-3\pi i/4}exp\left(3\pi i(\frac{a+d}{12c}-s(d,c)\right).$
\vspace{0.2cm}
\\
We adopt a similar approach to the method used in Theorem 1 to prove theorem 2. We only show the transformation law of $\theta_3$. As for $\theta_2$, and $\theta_4$  the proof is the same.\\
The generators for $\Gamma(2)$ are $T^2=\left(\begin{array}{cc}
   1  & 2 \\
   0  &1 \\
\end{array}\right)$ and $S=\left(\begin{array}{cc}
   1  & 0 \\
   2  &1 \\
\end{array}\right)$.
\vspace{0.25cm}\\
We now present three lemmas that are needed to prove Theorem 2.
\vspace{0.25cm}
\\
\textbf{Lemma 3.} $\epsilon_1(AT^{m})=\epsilon(A).e^{\pi im/2}$\\
\textit{Proof. } We will be using the matrix $AT^{2m}=\left(\begin{array}{cc}
   a  & 2am+b \\
   c  &2cm+d \\
\end{array}\right).$  \[\epsilon_1(AT^{2m})=e^{-3\pi i/4}exp\left(3\pi i \left( \frac{a+2cm+d}{12c}-s(d,c)\right)\right)~~~~~~~~~~~~~~~~~~~~~~~~~~~~~~~~~~~~~~~~~~~~~\]\[~~~=e^{-3\pi i/4}exp\left(3\pi i \left( \frac{a+d}{12c}-s(d,c)\right)\right).e^{\pi im/2}~~~~~~~~~~~~~~~~~~~~\]\[=\epsilon_1(A).e^{\pi i m/2}.~~~~~~~~~~~~~~~~~~~~~~~~~~~~~~~~~~~~~~~~~~~~~~~~~~~~~~~~\]
\vspace{0.35cm}
\\
\textbf{Lemma 4.} $\epsilon_1(AS)=\left\lbrace \begin{array}{ll}
\epsilon(A).e^{-\pi i/2}~~~~if~c+2d>0\\
\epsilon(A)~~~~if~c+2d<0\\
\end{array}\right.$
\vspace{0.3cm}
\\
\textit{Proof.} We have $AS=\left(\begin{array}{cc}
   a+2b  & b \\
   c+2d  &d \\
\end{array}\right)$. First we will treat the case $c+2d>0$,
\vspace{0.3cm}
\\
\[\epsilon_1(AS)=e^{-3\pi i/4}exp\left(3\pi i \left( \frac{a+2b+d}{12(c+2d)}-s(d,c+2d)\right)\right).\] Note also that if $d<0$, it will degenerate to the same result if one chooses $-s(d,c+2d)=s(-d,c+2d)$ so that when we use the reciprocity law we ensure that $-d>0$.\\
Using reciprocity law and the property of Dedekind sum, we have
\[s(d,c+2d)+\underbrace{s(c+2d,d)}_\textrm{=s(c,d)}=\frac{d}{12(c+2d)}+\frac{c+2d}{12d}-\frac{1}{4}+\frac{(a+2b)d-b(c+2d)}{12d(c+2d)}.\]
As a result
\[ \frac{a+2b+d}{12(c+2d)}-s(d,c+2d)=s(c,d)-\frac{c+2d}{12d}+\frac{1}{4} +\frac{b}{12d}.\]
Using reciprocity law again, we have
\[s(d,c)+s(c,d)=\frac{c}{12d}=\frac{d}{12c}-\frac{1}{4} +\frac{ad-bc}{12cd}.\]
We get
\[ s(d,c)-\frac{c+2d}{12d}+\frac{1}{4}+\frac{b}{12d}=-s(d,c)+\frac{a+d}{12c}-\frac{1}{6}.\]
Hence \[\epsilon_1(AS)=e^{-3\pi i/4}exp\left(3\pi i \left(\frac{a+d}{12c}-s(d,c)\right)\right).e^{-\pi i/2}\]\[=\epsilon_1(A).e^{-\pi i/2}.~~~~~~~~~~~~~~~~~~~~~~~~~~~~~~\]
For the case when $c+2d<0$, we take the matrix $AS=\left(\begin{array}{cc}
  -a-2b  & -b \\
   -c-2d  &-d \\
\end{array}\right)$ and the proof follows the same way where one uses the reciprocity law of Dedekind sum twice.
\vspace{0.25cm}
\\
\textbf{Lemma 5.}  \[\theta_3\left(\frac{z}{2\tau+1},\frac{\tau}{2\tau+1}\right)=i.e^{-\pi i /2}(2\tau+1)^{1/2}e^{\frac{2\pi i z^2}{2\tau+1}}\theta_3\left(z,\tau \right).\]
We have 
\[~~~~~~~~~~\theta_3(\frac{z}{\tau},\frac{-1}{\tau})=(-i\tau)^{1/2}e^{\frac{\pi i z^2}{\tau}}\theta_3(z,\tau).~~~~~~~~~(*)\]
Doing a change of variable $\tau\to \frac{-(2\tau+1)}{\tau}$ we get:
\[\theta_3(\frac{-z\tau}{2\tau+1},\frac{\tau}{2\tau+1})=\left(i.\frac{2\tau+1}{\tau}\right)^{1/2}e^{\frac{-\pi i z^2\tau}{2\tau+1}}\theta_3\left(z,-\frac{(2\tau+1)}{\tau}\right).\]
Another change of variable $z\to-z/\tau$, we obtain
\[~~~~~~~~~\theta_3(\frac{z}{2\tau+1},\frac{\tau}{2\tau+1})=\left(i.\frac{2\tau+1}{\tau}\right)^{1/2}e^{\frac{-\pi i z^2}{\tau(2\tau+1)}}\theta_3\left(\frac{-z}{\tau},-2-\frac{1}{\tau}\right).~~~~~~~~~~~~~~(**)\]
However $\theta_3(-z,\tau)=\theta_3(z,\tau)$ and $\theta_3(z,\tau+2m)=\theta_3(z,\tau)$  (see [3-X]). Hence ($**$) becomes:
\[\theta_3(\frac{z}{2\tau+1},\frac{\tau}{2\tau+1})=\left(i.\frac{2\tau+1}{\tau}\right)^{1/2}e^{\frac{-\pi i z^2}{\tau(2\tau+1)}}\theta_3\left(\frac{z}{\tau},-\frac{1}{\tau}\right).\]
Using (*)
\[~~~~~~~~~=\left(i.\frac{2\tau+1}{\tau}\right)^{1/2}e^{\frac{-\pi i z^2}{\tau(2\tau+1)}}(-i\tau)^{1/2}e^{\frac{\pi i z^2}{\tau}}\theta_3(z,\tau)\]
\[~~=(2\tau+1)^{1/2}e^{\frac{2\pi i z^2}{2\tau+1}}\theta_3(z,\tau)~~~~~~~~~~~~~~~~~~\]
\[=i.e^{-\pi i /2}(2\tau+1)^{1/2}e^{\frac{2\pi i z^2}{2\tau+1}}\theta_3(z,\tau).~~~~\]
which completes our proof of Lemma 5.
\vspace{0.4cm}
\\
We now present the proof for Theorem 2. Recall that we only have the theorem for $\theta_3$ since the proof follows in the same way for $\theta_2$ and $\theta_4$.\\

\textit{Proof of Theorem 2}  $AT^{2m}$ we replace $\tau\to T^{2m}\tau=\tau+2m$ in (8), we get
\[\theta_3\left(\frac{z}{c\tau+2cm+d},\frac{a\tau+2am+b}{c\tau+2cm+d}\right)~~~~~~~~~~~~~~~~~~~~~~~~~~~~~~~~~~~~~~~~~~~~~~\]
\[~~~~~~~~~~~~~~~~~~~~~~~~~~~~~~~~~~~~=i^{(d-1)(c/2+1)-\frac{b}{2}(a)+c/2}\epsilon_1(A)(c\tau+2cm+d)^{1/2}e^{\frac{\pi i z^2 c}{c\tau+2cm+d}}\theta_3(z,\tau+2m).~~~~~~(***)\]
Using the fact that $\theta_3(z,\tau+2m)=\theta_3(z,\tau)$ (see [3-X]), we have that (***) is equivalent to :
\[\theta_3\left(\frac{z}{c\tau+2cm+d},\frac{a\tau+2am+b}{c\tau+2cm+d}\right)=i^{(d-1)(c/2+1)-\frac{b}{2}(a)+c/2}\epsilon_1(A)(c\tau+2cm+d)^{1/2}e^{\frac{\pi i z^2 c}{c\tau+2cm+d}}\theta_3(z,\tau).~~~~~~~~~(9)\]
Now using the matrix $AT^{2m}=\left(\begin{array}{cc}
   a  & 2am+b \\
   c  &2cm+d \\
\end{array}\right)$, we must show that (9) is equivalent to
\[\theta_3\left(\frac{z}{c\tau+2cm+d},\frac{a\tau+2am+b}{c\tau+2cm+d}\right)~~~~~~~~~~~~~~~~~~~~~~~~~~~~~~~~~~~~~~~~~~~~~~~~~~~~~~~~~\]
\[~~~~~~~~~~~~~~~~~~~~~~~~~~~~~~~~~~~~=i^{(2cm+d-1)(c/2+1)-\frac{2am+b}{2}(a)+c/2}\epsilon_1(AT^{2m})(c\tau+2cm+d)^{1/2}e^{\frac{\pi i z^2 c}{c\tau+2cm+d}}\theta_3(z,\tau).\]
We are required to prove
\[~~~~~~~~~i^{(d-1)(c/2+1)-\frac{b}{2}(a)+c/2}=e^{m\pi i/2}.i^{(2cm+d-1)(c/2+1)-\frac{2am+b}{2}(a)+c/2}.~~~~~~~~~(10)\]
Note that \[\displaystyle e^{m\pi i/2}.i^{(2cm+d-1)(c/2+1)-\frac{2am+b}{2}(a)+c/2}~~~~~~~~~~~~~~~~~~~~~~~~\]\[~~~~~~~~~~~~~~~~~~~~~~~~~~~=i^{m+(2cm+d-1)(c/2+1)-\frac{2am+b}{2}(a)+c/2}\]\[~~~~~~~~~~~~~~~~~~~~~~~~~~~~~~~~~~~~~=i^{m+2cm(c/2+1)-a^2m}.i^{(d-1)(c/2+1)-\frac{b}{2}(a)+c/2}.\] We have to prove that \[\displaystyle i^{m+2cm(c/2+1)-a^2m}=1\]
\vspace{0.3cm}
\\
which is obvious since \[i^{m+2cm(c/2+1)-a^2m}=i^{m((c+1)^2-a^2)}.\] 
We have $c\equiv~ 0~(mod~ 2)$ which implies $(c+1)^2\equiv 1 ~(mod ~4)$. Similarly $-a^2\equiv -1~(mod~4)$, and thus
$(c+1)^2-a^2\equiv 0~(mod~4)$.\\ Hence \[i^{m((c+1)^2-a^2)}=i^{m+2cm(c/2+1)-a^2m}=1\].\\
This proves (10). Using Lemma 3, \[e^{m\pi i/2}\epsilon_1(A)=\epsilon_1(AS),\] we obtain that (9) is equivalent to applying the matrix $AT^{2m}=\left(\begin{array}{cc}
   a  & 2am+b \\
   c  &2cm+d \\
\end{array}\right)$ to $(8)$.
\vspace{0.35cm}
\\
For the case of $AS=\left(\begin{array}{cc}
   a+2b  & b \\
   c+2d  &d \\
\end{array}\right)$, assuming $c+2d>0$, we introduce the notation $\alpha(A)=i^{(d-1)(c/2+1)-\frac{b}{2}(a)+c/2}$ and hence we are required to prove that
\[\theta_3\left(\frac{z}{cS\tau+d},\frac{aS\tau+b}{c\tau+d}\right)=\alpha(A)\epsilon_1(A)(cS\tau+d)^{1/2}e^{\frac{\pi i z^2 c}{cS\tau+d}}\theta_3(z,S\tau)~~~~~~~~~(11)\]
is equivalent to
\[\theta_3\left(\frac{z}{(c+2d)\tau+d},\frac{(a+2b)\tau+b}{(c+2d)\tau+d}\right)=\alpha(AS)\epsilon_1(AS)((c+2d)\tau+d)^{1/2}e^{\frac{\pi i z^2 (c+2d)}{(c+2d)\tau+d}}\theta_3(z,\tau)~~~~~~~(12)\]
where $\alpha(AS)=i^{(d-1)(\frac{c+2d}{2}+1)-\frac{b}{2}(a+2b)+\frac{c+2d}{2}}.$
\vspace{0.25cm}
\\
Now \[~~~~~~~~~~~~~~~~~~~~~~~~~~~~~~~~\theta_3\left(\frac{z}{cS\tau+d},\frac{aS\tau+b}{c\tau+d}\right)=\theta_3\left(\frac{z(2\tau+1)}{(c+2d)\tau+d},\frac{(a+2b)\tau+b}{(c+2d)\tau+d}\right).~~~~~~~~~~~~~~~~~~~~(13)\]\\
We make the change of variable $\displaystyle z\to \frac{z}{2\tau+1}$, and (13) becomes \[\theta_3\left(\frac{z}{(c+2d)\tau+d},\frac{(a+2b)\tau+b}{(c+2d)\tau+d}\right)=\alpha(A)\epsilon_1(A)\left(\frac{(c+2d)\tau+d}{2\tau+1}\right)^{1/2}e^{\frac{\pi i z^2 c}{(2\tau+1)((c+2d)\tau+d)}}\theta_3(\frac{z}{2\tau+1},\frac{\tau}{2\tau+1}).\]
We now use Lemma 5
\[~~~~~~~~~~~~~~~~~~~~~~~~~~~~\theta_3\left(\frac{z}{2\tau+1},\frac{\tau}{2\tau+1}\right)=i.e^{-\pi i /2}(2\tau+1)^{1/2}e^{\frac{2\pi i z^2}{2\tau+1}}\theta_3\left(z,\tau \right).~~~~~~~~~~~~~~~~~~~~~(14)\]
This reduces the above into
\[\theta_3\left(\frac{z}{(c+2d)\tau+d},\frac{(a+2b)\tau+b}{(c+2d)\tau+d}\right)~~~~~~~~~~~~~~~~~~~~~~~~~~~~~~~~~~~~~~~~~~~~~~~~~~~~~~~~~~~\]
\[~~~~~~~~~~~~~~~~~~~~~=i^{(d-1)(c/2+1)-\frac{b}{2}(a)+c/2}.i.e^{-\pi i/2}\epsilon_1(A)\left(c+2d)\tau+d\right)^{1/2}e^{\frac{\pi i z^2 c}{(2\tau+1)((c+2d)\tau+d)}}e^{\frac{2\pi i z^2}{2\tau+1}}\theta_3(z,\tau)\]
\[~~~~~~~~~~~~~~~~~~~~~=i^{(d-1)(c/2+1)-\frac{b}{2}(a)+c/2+1}\epsilon_1(AS)\left((c+2d)\tau+d\right)^{1/2}e^{\frac{\pi i z^2(c+2d)}{(c+2d)\tau+d}}\theta_3(z,\tau).\]
Comparing to (12), we just have to prove that
\[~~~~~~~~~~~~~~~\alpha(AS)=i^{(d-1)(\frac{c+2d}{2}+1)-\frac{b}{2}(a+2b)+\frac{c+2d}{2}}=i^{(d-1)(c/2+1)-\frac{b}{2}(a)+c/2+1}.~~~~~~~~~~~~~~~~~~~~~~~(15)\]
Expanding and collecting terms, we get
\[~~~~~~~~~~~~~~~~~~~~~~~~~~~~~~~i^{(d-1)(\frac{c+2d}{2}+1)-\frac{b}{2}(a+2b)+\frac{c+2d}{2}}=i^{(d-1)(c/2+1)-\frac{b}{2}(a)+c/2}.i^{(d^2-b^2)}.~~~~~~~~~~~~~~~(16)\]
But $d^2\equiv 1~(mod~4)$ and $b^2\equiv 0~(mod~4)$, so $i^{(d^2-b^2)}=i$. This proves that (16) implies (15), which in turn proves that (11) and (12) are equivalent.\\
For the case of $c+2<d$ the approach would be the same as in Theorem 1. Using $AS=\left(\begin{array}{cc}
   -a-2b  &-b \\
  -c-2d  &-d \\
\end{array}\right)$, where one makes use of the fact that $\theta_3$ is even in terms of $z$, this completes the proof of Theorem 2.

\vspace{1cm}
\begin{center}
    ACKNOWLEDGEMENT
\end{center}
We would like to express our deepest gratitude to our advisor Professor Wissam Raji. We also like to thank the Center of Advanced Mathematical Sciences (CAMS) at the American university of Beirut (AUB) for the guidance and support we are recieving from the summer research camp (SRC).
\\

Department of Mathematics, American University of Beirut, Beirut, Lebanon
\\
\textit{E-mail address: mmm133@mail.aub.edu}\\
\textit{E-mail address: ays11@mail.aub.edu}

\end{document}